%% file: Burns--Krantz_Rigidity_arXiv_.tex
\documentclass[11pt,letterpaper,reqno]{amsart}
\usepackage[T1]{fontenc}
\usepackage{lmodern}
\usepackage[margin=1.1 in]{geometry}
\usepackage{amsthm}
\usepackage{amsfonts,amssymb}
\usepackage{comment}
\usepackage{enumitem,chngcntr}
\usepackage{xcolor}

\usepackage{hyperref}
\newtheorem{theorem}{Theorem}[subsection]
\counterwithin*{theorem}{section}   
\renewcommand{\thetheorem}{
  \ifnum\value{subsection}=0 \thesection\else\thesubsection\fi.\arabic{theorem}}
\newtheorem{lemma}[theorem]{Lemma}
\newtheorem{proposition}[theorem]{Proposition}

\theoremstyle{definition}
\newtheorem{definition}[theorem]{Definition}
\newtheorem{example}[theorem]{Example}

\theoremstyle{remark}
\newtheorem{remark}[theorem]{Remark}
\newcommand{\BB}{{\mathbb B}}

\newcommand{\DD}{{\mathbb D}}
\newcommand{\OO}{{\mathcal O}}

\newcommand{\NN}{{\mathbb N}}

\newcommand{\RR}{{\mathbb R}}
\newcommand{\CC}{{\mathbb C}}
\newcommand{\cC}{{\mathcal C}}

\newcommand{\ole}{\partial \mathbb{D}}
\newcommand{\co}{\overline}

\newcommand{\cce}{\mathcal{E}(p)}

\DeclareMathOperator{\Aut}{Aut} \DeclareMathOperator{\id}{id}

 \DeclareMathOperator{\dist}{dist}

\renewcommand{\phi}{\varphi}
\title{Burns--Krantz rigidity for convex complex ellipsoids}
\author{Anand Chavan}\email{anand.chavan@doctoral.uj.edu.pl; anandmaths57@gmail.com}
\address{Doctoral School of Exact and Natural Sciences, Institute of Mathematics, Faculty of Mathematics and Computer Science, Jagiellonian
University,  \L ojasiewicza 6, 30-348 Krak\'ow, Poland}

\begin{document}
\begin{abstract}
      In this article we study the Burns--Krantz rigidity phenomenon. We prove Burns--Krantz rigidity at every boundary point of the convex complex ellipsoids $\cce\subset \CC^n$ with $p_j\in\{1/2,1\}$ for all $j$. This extends the result of \cite{Gen-Mig 1997}, where, for $p_j=1/2$ for all $j$, rigidity was established at every boundary point when $n=2$ (the diamond), but only at smooth boundary points when $n\geq 3$. Our proof is based on the approach of \cite{Zwo 2025a}. 
\end{abstract}
\maketitle

\section{Introduction and Background}
Recall the following classical result of Cartan \cite{Cartan1931}:
\begin{theorem}[Cartan's Uniqueness Theorem, 1931]
   Let $D\subset \CC^n$ be a bounded domain, and let $f:D\to D$ be a holomorphic map. If there exists $z_0\in D$ such that \[f(z)=z+o(\|z-z_0\|) \quad \text{as }z\to z_0,\] then $f\equiv \id_{D}$.
\end{theorem}
 A natural question is whether an analogous rigidity phenomenon is observed when the fixed point lies on the boundary of the domain rather than in its interior. More precisely, can one conclude that a holomorphic self-map is the identity if it agrees with the identity map to sufficiently high order at a boundary point? In 1994, Burns and Krantz answered this question affirmatively for smoothly bounded strongly pseudoconvex domains, provided that the first-order contact condition is replaced by a higher-order one. We recall their result below.

\begin{theorem}[\cite{Bur-Kra 1994}]
Let $D \subset \CC^n$ be a bounded strongly pseudoconvex domain with a $\cC^{6}$-smooth boundary, and let $b \in \partial D$. Suppose $f \colon D \to D$ is holomorphic and satisfies
\[
f(z)=z+o\!\left(\|z-b\|^3\right)\qquad \text{as } z\to b.
\]   
Then $f \equiv \id_{D}$.
\end{theorem}
This result was new even in $\CC$ for the unit disc. The exponent $3$ is optimal, as an example in the case of the unit disc shows (see \cite[Remark~1]{Bur-Kra 1994}). Recently, Zimmer \cite[Theorem 4.1]{Zim 2022} gave another proof for the unit disc using new estimates for the Poincaré metric \cite{Chr-Sho 2019}.

Huang proved the following, in which ``local'' strong pseudoconvexity is used instead of global strong pseudoconvexity.

\begin{theorem}[{\cite[Theorem 2.5]{Hua 1995}}]
    Let $D$ be a bounded domain in $\CC^n$. Let $b\in \partial D$ be a strongly pseudoconvex point. Let $f:D\to D$ be a holomorphic map such that
\[
f(z)=z+o\!\left(\|z-b\|^3\right)\qquad \text{as } z\to b.
\]
Then $f \equiv \id_{D}$.
\end{theorem}

After the work of Burns and Krantz \cite{Bur-Kra 1994}, several authors have generalized the result to various settings, for instance by assuming different regularity conditions either on the map or on the boundary of the domain. We refer the reader to the recent articles \cite{Ron 2026,Zim 2022} for a detailed account of the developments in this direction. Let us now recall some important results in our context.

\begin{theorem}[\cite{Zim 2022}]\label{thm:zimbk}
    Let $D\subset\CC^n$ be a bounded convex domain with a $\cC^2$-smooth boundary, and let $f:D\to D$ be a holomorphic map. Suppose that for some $b\in \partial D$, \[f(z)=z+o(\|z-b\|^4)\quad \text{as } z\to b.\] Then $f\equiv \id_{D}$.
\end{theorem}

While this article was being completed, Lou, Qiao, Wang and
Zhong~\cite{Lou-Qia-Wan-Zho 2026} showed that the exponent~$4$ is not optimal,
answering Zimmer's question: they proved Burns--Krantz rigidity with exponent~$3$
for bounded convex domains with $\cC^{1,\alpha}$-smooth boundary.

Before we move further, let us formalize some notions. 
\begin{definition}[Burns--Krantz map]
    Let $D$ be a domain in $\CC^n$, $b \in \partial{D}$, and $\gamma\geq 3$. If a holomorphic mapping $f:D \to D$ satisfies \[f(z)= z + o\bigl( \|z-b\|^\gamma\bigr) \qquad \text{as } z\to b,\] then we say that $f$ is a \emph{Burns--Krantz map} at the boundary point $b$ with exponent $\gamma$.
\end{definition}

\begin{definition}[Burns--Krantz rigidity property]
    Let $D$ be a domain in $\CC^n$, $b\in\partial D$ and $\gamma\geq 3$. We say that $(D,b)$ has \emph{the Burns--Krantz rigidity property} with exponent $\gamma$ if every Burns--Krantz map $f:D\to D$ at $b$ with exponent $\gamma$ is the identity map.
\end{definition}

Throughout, we consider Burns--Krantz maps with exponent~$3$, unless otherwise stated.

Most results on Burns--Krantz rigidity are highly sensitive to the boundary regularity of the domain. In \cite{Zwo 2025a}, Lempert's techniques are used to understand the Burns--Krantz rigidity phenomenon. In fact, that approach works with a complex geodesic and its holomorphic left inverse, i.e.\ with the internal geometry of the domain. The following lemma from \cite{Zwo 2025a} is the key tool in this approach.
\begin{lemma}[Zwonek \cite{Zwo 2025a}]\label{lem:zwobk}
Let $D$ be a domain in $\CC^n$, $b \in \partial D$, and let 
$f: D \to D$ be a Burns--Krantz map at $b$. Assume that $\phi: \DD \to D$ is a complex geodesic that is Lipschitz continuous at $1$ with $\phi(1) = b$, and let $G: D \to \DD$ be a left inverse of $\phi$ that is Lipschitz continuous at $b$. Then $G$ is a left inverse of $f \circ \phi$. In particular, $f \circ \phi$ is a complex geodesic.
\end{lemma}
By combining this lemma with a detailed understanding of complex geodesics and their left inverses, Burns--Krantz rigidity was established at nonsmooth boundary points of $\DD^n$ and at points of the Shilov boundary of the symmetrized bidisc \cite{Zwo 2025a}. For smooth boundary points of $\DD^n$, the corresponding rigidity result had already been obtained by Forn\ae ss and Rong \cite{For-Ron 2021}. Together, these results establish Burns--Krantz rigidity at every boundary point of the polydisc.

Recall that the \emph{convex complex ellipsoid} $\cce\subset\CC^n$ for $p=(p_1,\dots,p_n)$ is defined by
\[\cce:=\left\{(z_1,\dots,z_n)\in \CC^n: \sum_{j=1}^{n}|z_{j}|^{2p_j}<1\right\},\qquad p_{j}\geq1/2 \text{ for all }j.\]
When $p_j=1/2$ for every $j$, the domain $\cce$ has a nonsmooth boundary. Its boundary is smooth precisely at points whose coordinates are all nonzero.

In \cite{Gen-Mig 1997}, Burns--Krantz rigidity was established at every boundary point (including nonsmooth ones) of $\cce\subset\CC^2$ when $p_1=p_2=1/2$. This domain is known as the \emph{diamond} and is denoted by $\triangle$. For $n\geq 3$, the result was established at the smooth boundary points of $\cce\subset\CC^n$ when $p_j=1/2$ for all $j$.

\begin{theorem}[\cite{Gen-Mig 1997}]\label{thm:gen_mig_bkr}
    Let $\cce\subset \CC^n$ be a convex complex ellipsoid with $p_j=1/2$ for all $j$, and let $b\in \partial \cce$.
    \begin{enumerate}
     \renewcommand{\labelenumi}{(\roman{enumi})}
     \item If $n=2$, then $(\cce,b)$ has the Burns--Krantz rigidity property. 
     \item If $n\geq 3$ and $b\in \partial\cce \cap (\CC\setminus\{0\})^n$, then $(\cce,b)$ has the Burns--Krantz rigidity property.
    \end{enumerate}    
\end{theorem} 
Huang \cite{Hua 1995} also studied the \emph{egg domains}, i.e.\ $\cce\subset \CC^2$ with $p_{1}=1$ and $p_2=k \in \NN$. In this case, however, all boundary points are smooth.  

The proof of Theorem~\ref{thm:gen_mig_bkr} given in \cite{Gen-Mig 1997} relies on explicit complex geodesics meeting the prescribed boundary point and on the construction of corresponding left inverses that extend holomorphically beyond the boundary, which is similar in essence to the approach in \cite{Zwo 2025a}.

Although the approaches in \cite{Zwo 2025a} and \cite{Gen-Mig 1997} rely on complex geodesics and their left inverses, obtaining a sufficiently precise description of these objects is itself a challenging problem.

In \cite{Jar-Pfl-Zei 1993} a complete description of the complex geodesics of the convex complex ellipsoids was given. The complex geodesics are unique (up to an automorphism of the unit disc) for any two distinct points of $\cce$. Moreover, when $p_{j}\in \{1/2,1\}$, complex geodesics extend holomorphically to a neighborhood of $\co{\DD}$. Following a similar approach, we establish Burns--Krantz rigidity at every boundary point of $\cce$ when $p_j\in\{1/2,1\}$ for all $j$. This extends the result of \cite{Gen-Mig 1997} to all boundary points in this class of domains and, when $p_j=1$ for every $j$, provides an essentially new proof of Burns--Krantz rigidity for the unit ball.

\begin{theorem}\label{thm:bkcce}
    Let $\cce\subset \CC^n$ be a convex complex ellipsoid with $p_{j}\in\{1/2,1\}$ for all $j\in \{1,\dots,n\}$, and let $b\in \partial \cce$. Then $(\cce,b)$ has the Burns--Krantz rigidity property.
\end{theorem}

 \begin{remark}
    As mentioned earlier, Lou, Qiao, Wang and Zhong \cite{Lou-Qia-Wan-Zho 2026} recently improved Zimmer's result (Theorem~\ref{thm:zimbk}). In particular, \cite[Theorem~1.1]{Lou-Qia-Wan-Zho 2026} yields Burns--Krantz rigidity with exponent~$3$ at every boundary point of $\cce$ whenever $p_j>1/2$ for all $j$, since in that case $\partial\cce$ is $\cC^{1,\alpha}$-smooth for some $\alpha\in(0,1]$. The case where $p_j=1/2$ for some $j$, however, is not covered by \cite[Theorem~1.1]{Lou-Qia-Wan-Zho 2026}, since at boundary points $b\in\partial\cce$ with $p_j=1/2$ and $b_j=0$ the boundary $\partial\cce$ is not even $\cC^{1}$ at $b$.

    Theorem~\ref{thm:bkcce} together with the result for the polydisc \cite{Zwo 2025a} and the diamond \cite{Gen-Mig 1997} gives examples of Burns--Krantz rigidity with exponent~$3$ at boundary points where the boundary regularity is lower than $\cC^{1,\alpha}$.

    \cite[Theorem~1.3]{Lou-Qia-Wan-Zho 2026} proves Burns--Krantz rigidity for bounded convex domains with an exponent $\gamma$ that is always strictly greater than~$3$. In particular, for $\cce$ with $p_j=1/2$ for some $j$ it does not give exponent~$3$.

    The proof of Theorem~\ref{thm:bkcce} presented here relies instead on the geometry of complex geodesics and the holomorphic extension of their Lempert left inverses. In particular, the results concerning the existence, extension, and uniqueness of these objects are of independent interest.
 \end{remark}
 \subsection{Sketch of the proof of Theorem~\ref{thm:bkcce}}\label{sketch}
\begin{enumerate}[ itemindent=0pt, topsep=2pt, itemsep=3pt, parsep=0pt, listparindent=0pt, leftmargin=1.6em]
    \item \textbf{Existence of complex geodesics.} Every complex geodesic of $\cce$ with $p_j\in\{1/2,1\}$ extends
    holomorphically to a neighborhood of $\co{\DD}$.

    Let $q$ be an arbitrary point of $\cce$. Then there exists a complex geodesic $\phi:\DD\to\cce$ such that $\phi(1)=b\in \partial \cce$ and $\phi(0)=q\in\cce$.

    \item \textbf{Existence of a holomorphic left inverse and application of
    Zwonek's lemma.}
    We will prove that there exists a left inverse $G:\cce\to\DD$ of $\phi$ such that
    $G\circ\phi=\id_{\DD}$ with $G(b)=1$ and $G$ is holomorphic near $b$.

    Then, applying Zwonek's lemma with the Burns--Krantz map $f$, we find that $\psi:=f\circ\phi$ is also a complex
    geodesic with the same left inverse $G$, and $\psi(1)=b$.

    \item \textbf{Equality of derivatives at the boundary.}
    We show that $\phi^{(l)}(1)=\psi^{(l)}(1)$ for $0\leq l\leq3$, and, using the explicit description of complex geodesics of $\cce$, we then show $\psi\equiv\phi$.
\end{enumerate}
From the above, we have $\psi(\lambda)=f(\phi(\lambda))=\phi(\lambda)$ for all $\lambda\in \DD$. In particular, $f(q)=f(\phi(0))=\phi(0)=q$. Since $q\in \cce$ was arbitrary, $f\equiv \id_{\cce}$.

The rest of the article is devoted to a detailed proof of Theorem~\ref{thm:bkcce}. In Section~\ref{sec:cg_li} we recall the necessary background on complex geodesics and Lempert left inverses. In Section~\ref{sec:proof} we expand each step of the above sketch into a separate subsection and then complete the proof. In Subsection~\ref{sec:exe_cg} we also study the (non)uniqueness of complex geodesics through two distinct points of $\co{\cce}$ with $p_j\in\{1/2,1\}$.

\section{Complex Geodesics and Left Inverses}\label{sec:cg_li}
In this section we recall the notions and results on complex geodesics and their left inverses that will be used in the proof of the main theorem. For background on invariant distances and metrics, complex geodesics and left inverses, we refer to \cite{Aba 1989,Jar-Pfl 2013}.

\begin{definition}
    Let $D\subset \CC^n$ be a domain. A holomorphic map $\phi:\DD\to D$ is called a \emph{complex geodesic} if there exists a holomorphic map $F:D\to \DD$, called a \emph{left inverse} of $\phi$, such that $F\circ\phi=\id_{\DD}$.
\end{definition}
Note that $\phi\circ a$ is again a complex geodesic for every $a\in\Aut(\DD)$, with left inverse $a^{-1}\circ F$. We say that a complex geodesic $\phi$ \emph{passes through} two points $z,w\in D$ if $z,w\in\phi(\DD)$, that is, $\phi(\lambda_1)=z$ and $\phi(\lambda_2)=w$ for some $\lambda_1,\lambda_2\in\DD$. Given $z\in D$ and $X\in\CC^n\setminus\{0\}$, we say that $\phi$ \emph{passes through $z$ in the direction of $X$} if, after composing with an automorphism of $\DD$ if necessary, we have $\phi(0)=z$ and $\phi'(0)=\alpha X$ for some $\alpha>0$. In this case $\phi$ is called an \emph{infinitesimal complex geodesic} for the datum $(z;X)\in D\times(\CC^n\setminus\{0\})$.

In \cite{Lem 1981,Lem 1982,Lem 1984}, Lempert proved that in every bounded convex domain, and also in bounded strongly linearly convex domains, any two distinct points lie on a complex geodesic. Moreover, the Carathéodory distance and metric coincide with the Kobayashi distance and metric, respectively. We refer to \cite{Kos-War 2013} for more details on Lempert's theorem.

We are also interested in the existence of complex geodesics through boundary points. More precisely, given two distinct points $z,w\in \co{D}$, is there a complex geodesic $\phi:\DD\to D$, extending continuously to $\co{\DD}$, with $z,w\in \phi(\co{\DD})$? In \cite{Cha-Hu-Lee 1988}, it was shown that for a bounded strongly linearly convex domain $D\subset \CC^n$, $n\geq2$, with $\cC^3$-smooth boundary, any two distinct points of $\co{D}$ lie on (the closure of) a complex geodesic, and Mercer \cite{Mer 1993} proved the same for $m$-convex domains, which we recall later. 

Similarly, one can pose the question about the existence of infinitesimal complex geodesics for the data $(z;X)\in \co{D}\times (\CC^n\setminus \{0\})$. For bounded convex or strongly linearly convex domains, existence for every $(z;X)\in D\times (\CC^n\setminus\{0\})$ follows from Lempert's theorem~\cite{Lem 1981,Lem 1982,Lem 1984}. When $z\in \partial D$ and $X\in \CC^n$ is an admissible direction, i.e.\ $X\in \CC^n\setminus T_{z}^{\CC}\partial D$, the existence is known for bounded strongly linearly convex domains with $\cC^3$-smooth boundary \cite{Hua 1994}, where $T_{z}^{\CC}\partial D$ denotes the complex tangent space at $z\in \partial D$.

Now we turn to uniqueness. Uniqueness is understood in the following sense:

\begin{definition}
Let $D\subset\CC^n$ be a bounded domain, and assume
that the complex geodesics under consideration admit
continuous extensions to $\co{\DD}$.

\begin{enumerate}
    \item \textbf{Uniqueness for two distinct points.}
    Let $z,w\in\co{D}$ be distinct.
    A complex geodesic $\phi:\DD\to D$ joining
    $z$ and $w$, that is,
    \[
        z,w\in\phi(\co{\DD}),
    \]
    is said to be \emph{unique} if, for every
    complex geodesic $\psi:\DD\to D$ satisfying
    \[
        z,w\in\psi(\co{\DD}),
    \]
    we have
    \[
        \phi(\DD)=\psi(\DD).
    \]

    \item \textbf{Uniqueness for a point and a direction.}
    Let $z\in\co{D}$ and
    $X\in\CC^n\setminus\{0\}$.
    Assume that the complex geodesics under consideration
    admit $\cC^1$-extensions at their boundary
    preimages whenever $z\in\partial D$.

    A complex geodesic through $z$ in the direction
    $X$ is said to be \emph{unique} if, for any two
    complex geodesics $\phi,\psi:\DD\to D$ satisfying,
    after suitable reparametrizations,
    \[
        \phi(\tau)=\psi(\tau)=z,
        \qquad
        \phi'(\tau),\psi'(\tau)\in A_X,
    \]
    where
    \[
        \begin{cases}
        \tau=0,\quad A_X=(\CC\setminus\{0\})X,
            & \text{if }z\in D,\\[2mm]
        \tau=1,\quad A_X=\RR_{>0}X,
            & \text{if }z\in\partial D,
        \end{cases}
    \]
    we have
    \[
        \phi(\DD)=\psi(\DD).
    \]
\end{enumerate}
\end{definition}
In both cases, $\phi(\DD)=\psi(\DD)$ if and only if $\phi=\psi\circ a$ for some $a\in\Aut(\DD)$. Indeed, if $F$ and $G$ are left inverses of $\phi$ and $\psi$ respectively, and $\phi(\DD)=\psi(\DD)$, then $a:=G\circ\phi$ and $F\circ\psi$ are holomorphic self-maps of $\DD$ that are inverse to each other, and $\psi\circ a=\phi$. Thus uniqueness of complex geodesics is always understood up to automorphisms of $\DD$.

In bounded strictly convex domains, as well as in bounded strongly linearly convex domains, any two distinct points of $D$ lie on a unique complex geodesic, and for every datum $(z;X)\in D\times (\CC^n\setminus\{0\})$ there exists a unique infinitesimal complex geodesic \cite{Lem 1982, Lem 1981,Lem 1984}. Moreover, if $D$ is a bounded strongly linearly convex domain with $\cC^3$-smooth boundary, then any two distinct points of $\co{D}$ are joined by a unique complex geodesic \cite{Cha-Hu-Lee 1988}. For infinitesimal complex geodesics, uniqueness also holds for every datum $(z;X)\in \partial D \times (\CC^n \setminus T_{z}^{\CC}\partial D)$ \cite{Hua-Wan 2022}.

The existence of (infinitesimal) complex geodesics holds true for every bounded convex domain, but uniqueness may fail without strict convexity. The bidisc $\DD^2$ is the simplest example. The convex complex ellipsoids form a notable class of convex, not necessarily strictly convex, domains in which complex geodesics are nevertheless unique for any two distinct interior points and interior infinitesimal data $(z;X)\in \cce\times (\CC^n\setminus \{0\})$. Jarnicki, Pflug and Zeinstra \cite{Jar-Pfl-Zei 1993} gave a complete description of the complex geodesics of convex complex ellipsoids, one of the few classes of convex domains for which such a description is known. The following theorem is a reformulation of their result.

\begin{theorem}[\cite{Jar-Pfl-Zei 1993}]\label{thm:geocce}
    Let $\cce\subset\CC^n$ be a convex complex ellipsoid, i.e.\ $p_j\geq 1/2$ for all $j$. A holomorphic map $\phi=(\phi_1,\dots,\phi_n):\DD\to \CC^n$ with $\phi_{j}\not \equiv 0$ for all $j$ is a complex geodesic of $\cce$ if and only if
    \begin{equation}\label{eq:CCE_geo}
        \phi_{j}(\lambda)=a_j\left(\frac{\lambda - \alpha_j}{1 - \overline{\alpha_j}\lambda}\right)^{r_j}
        \left(\frac{1 - \overline{\alpha_j}\lambda}{1 - \overline{\alpha_0}\lambda}\right)^{1/p_j},
        \qquad j=1,\dots,n,
    \end{equation}
    where
    \begin{enumerate}
        \item $a_1,\dots,a_n\in\CC\setminus\{0\}$, $\alpha_0\in\DD$, and $r_j\in\{0,1\}$ for $j=1,\dots,n$;
        \item $\alpha_j\in\DD$ if $r_j=1$, and $\alpha_j\in\co\DD$ if $r_j=0$;
        \item the constants satisfy
        \[
            \alpha_0 = \sum_{j=1}^{n}|a_j|^{2p_j}\alpha_{j},
            \qquad
            1+|\alpha_{0}|^2 = \sum_{j=1}^{n}|a_j|^{2p_j}(1+|\alpha_j|^2);
        \]
        \item the case $r_1=\dots=r_n=0$, $\alpha_0=\alpha_1=\dots=\alpha_n$ is excluded;
        \item the branches of the powers are chosen so that $1^{1/p_j}=1$, $j=1,\dots,n$.
    \end{enumerate}
\end{theorem}
The assumption that $\phi_j\not\equiv 0$ is not restrictive. If $\phi_j\equiv 0$ for some $j$, then the mapping $\hat{\phi}=(\phi_{1},\dots,\phi_{j-1},\phi_{j+1},\dots,\phi_{n})$ is a complex geodesic in the lower-dimensional convex complex ellipsoid $\mathcal{E}(p')$, where $p'=(p_1,\dots,p_{j-1},p_{j+1},\dots,p_n)$. Hence no generality is lost by considering geodesics with no component identically zero. 

We will see in a later section that any two distinct points of $\co{\cce}$ lie on (the closure of) a complex geodesic. Moreover, when $p_{j}\in \{1/2,1\}$, every complex geodesic of $\cce$ extends holomorphically to a neighborhood of $\co{\DD}$. Surprisingly, despite this regularity, complex geodesics through boundary points need not be unique, although any two distinct points of $\cce$ lie on a unique complex geodesic.

We now turn to left inverses. In his proof that the Carathéodory distance and metric coincide with the Kobayashi distance and metric, respectively, on bounded strongly linearly convex domains, Lempert \cite{Lem 1984} implicitly constructed a left inverse of the complex geodesic passing through any two distinct points. Royden and Wong \cite{Roy-Won 1983} gave another proof of Lempert's theorem for bounded convex domains (see also \cite[Section~2.6]{Aba 1989}). Their proof also constructs left inverses of complex geodesics by means of dual maps. We briefly describe this construction of the Lempert left inverse of a complex geodesic from its dual map on a bounded convex domain.

Let $D\subset\CC^n$ be a bounded convex domain, let $\phi:\DD\to D$ be a complex geodesic, and let $\widetilde{\phi}$ be a dual map of $\phi$. For each $z\in D$ the equation
\begin{equation}
    \langle z-\phi(\lambda),\widetilde{\phi}(\lambda)\rangle=0
\end{equation}
where \[
    \langle a,b\rangle=\sum_{j=1}^n a_jb_j,
    \qquad a,b\in\CC^n,
\]
has a unique solution $\lambda=\lambda(z)\in\DD$, and $F(z):=\lambda(z)$ defines a holomorphic map $F:D\to\DD$ such that
\begin{equation}\label{eq:dual_lem_inv}
     \langle z-\phi(F(z)),\widetilde{\phi}(F(z))\rangle=0 \quad\text{for all } z\in D,
\end{equation}
in particular $F\circ\phi=\id_{\DD}$ (see \cite[Proposition~2.6.22]{Aba 1989}). By construction, every fiber of $F$ is the intersection of $D$ with an affine complex hyperplane:
\[F^{-1}(\lambda)\cap D= H_\lambda \cap D, \quad \lambda\in\DD\]
where
\[
    H_\lambda=\left\{z\in\CC^n:\langle z-\phi(\lambda),\widetilde{\phi}(\lambda)\rangle=0\right\},
    \quad \lambda\in\DD.
\]
A complex geodesic may have left inverses whose fibers are not affine.

\begin{example}\label{ex:leftinv_dia}
    Consider the diamond $\triangle=\{|z_1|+|z_2|<1\}\subset \CC^2$, that is, $\cce$ with $p=(1/2,1/2)$, and the complex geodesic $\phi(\lambda)=(\lambda/2,1/2)$. The map $F(z_1,z_2)=z_1/(1-z_2)$ is a left inverse of $\phi$: since $|z_1|<1-|z_2|\leq|1-z_2|$ on $\triangle$, $F$ maps $\triangle$ into $\DD$, $F\circ\phi=\id_{\DD}$, and $F^{-1}(\lambda)=\{z_1=\lambda(1-z_2)\}\cap\triangle$ is affine. 
    
    The map $\widehat{F}(z_1,z_2)=4z_1z_2$ is also a left inverse of $\phi$, since $|4z_1z_2|\leq(|z_1|+|z_2|)^2<1$ on $\triangle$ and $\widehat F\circ\phi=\id_{\DD}$, but its fibers $\{z_1z_2=\lambda/4\}$ are not affine.
\end{example}
Among all left inverses of a complex geodesic, those with affine fibers are called Lempert left inverses. We define them formally below.
\begin{definition}
    Let $D\subset \CC^n$ be a domain and let $\phi:\DD\to D$ be a complex geodesic. A left inverse $F$ of $\phi$ is called a \emph{Lempert left inverse} if its fibers are affine, that is, for every $\lambda\in\DD$ there is an affine complex hyperplane $H_\lambda\subset\CC^n$ such that
    \begin{equation}\label{eq:lem_inv_hp}
        F^{-1}(\lambda)\cap D= H_\lambda \cap D.
    \end{equation}
\end{definition}

\section{Proof}\label{sec:proof}
In this section we develop the proof of the main theorem, Theorem~\ref{thm:bkcce}. Each subsection represents a step in the sketch of the proof described earlier in Subsection~\ref{sketch}.
\subsection{Existence of a complex geodesic}\label{sec:exe_cg}
Theorem~\ref{thm:geocce} gives a complete description of the complex geodesics of $\cce$. In fact, they are unique (up to an automorphism of the unit disc) for any two distinct points of $\cce$. A natural question is whether each boundary point can be joined to some interior point by a complex geodesic, even when one assumes that all complex geodesics admit holomorphic extensions beyond $\co{\DD}$.

In the symmetrized bidisc, every complex geodesic does extend holomorphically to a neighborhood of $\co\DD$, yet complex geodesics fail to connect any topological boundary points other than those lying in the Shilov boundary \cite{Pfl-Zwo 2005}.

Mercer \cite{Mer 1993} answered this question for $m$-convex domains.

\begin{definition}[$m$-convex domain]\label{def:m_cvx}
    Let $D\subset \CC^n$ be a bounded convex domain. For
    $z\in D$ and $v\in\CC^n\setminus\{0\}$, put
    \[
    \delta_{D}(z):=\dist(z,\partial D)
    \quad\text{and}\quad
    \delta_{D}(z;v)
    :=\sup\Bigl\{r>0 : z+\lambda \tfrac{v}{\|v\|}\in D \text{ for all } |\lambda|<r\Bigr\}.
    \]
    We say that $D$ is \emph{$m$-convex} if there exists a constant $C>0$
    such that
    \[
    \delta_{D}(z;v)\leq C\,\delta_{D}(z)^{1/m}
    \]
    for every $z\in D$ and every $v\in\CC^n\setminus\{0\}$.
\end{definition}

\begin{theorem}[\cite{Mer 1993}]\label{thm:mer_93}
Let $D\subset\CC^n$ be a bounded $m$-convex domain. For any two distinct points $z,w\in \co D$ there is a complex geodesic whose continuous extension to $\co\DD$ contains these two points.
\end{theorem}

Next we observe that $\cce$ is an $m$-convex domain. We defer the proof to Appendix~\ref{Appendix}.

\begin{proposition}\label{prop:cce_m_cnv}
    $\cce$ is $m$-convex with $m=\max\{2,2p_1,\dots,2p_n\}$.
\end{proposition}

As a consequence, any two distinct points in $\co{\cce}$ can be joined by a complex geodesic. In particular, for any $b\in \partial \cce$ and $q\in \cce$ there exists a complex geodesic $\phi:\DD\to \cce$, extending continuously to $\co\DD$, such that, up to an automorphism of the unit disc, $\phi(1)=b$ and $\phi(0)=q$. Moreover, when $p_j\in\{1/2,1\}$, the complex geodesic extends holomorphically to a neighborhood of $\co{\DD}$.

In the diamond, any two distinct points of the closure lie on a complex geodesic that extends holomorphically to a neighborhood of $\co{\DD}$. By \cite{Jar-Pfl-Zei 1993}, this geodesic is unique when both points lie in the interior. Once boundary points are allowed, however, uniqueness may fail, as the following two examples show.

\begin{example}\label{eg:cce_cg_nu} For the diamond there exist complex geodesics with $\phi(i)=\psi(i)$, $\phi(-i)=\psi(-i)$ and $\phi(\DD)\neq\psi(\DD)$.

    Consider
    \[
\phi_t(\lambda)
=
\left(
\frac{(\lambda+t)(1+t\lambda)}{2(1+t^2)},
\frac{(\lambda-t)(1-t\lambda)}{2(1+t^2)}
\right),
\qquad t\in(-1,1).
\]
 Observe that for every $t\in (-1,1)$, $\phi_{t}$ is a complex geodesic (Theorem~\ref{thm:geocce}) with $F(z_1,z_2)=z_1+z_2$ as a left inverse, $F\circ \phi_{t}=\id_{\DD}$ and 
 \[\phi_{t}(i)=\left(\frac{i}{2},\frac{i}{2}\right)=:w_+,\qquad \phi_{t}(-i)=\left(-\frac{i}{2},-\frac{i}{2}\right)=:w_-.\]
 So for every $t\in (-1,1)$, $w_+,w_-\in \phi_{t}(\co{\DD})$.

 We will show that $\phi_{t_1}(\DD)\neq\phi_{t_2}(\DD)$ whenever $t_1\neq t_2$.

 Suppose, on the contrary, that $\phi_{t_{1}}(\DD)=\phi_{t_2}(\DD)$ for some $t_1\neq t_2$. Then there exists $\sigma\in \DD$ such that \[\phi_{t_1}(0)=\phi_{t_2}(\sigma).\]
 Since $F$ is a left inverse of both $\phi_{t_{1}}$ and $\phi_{t_2}$,
 \[0=F(\phi_{t_1}(0))=F(\phi_{t_2}(\sigma))=\sigma,\]
 therefore $\sigma=0$, and the equality of these two geodesics at $0$ gives \[\frac{t_1}{2(1+t_{1}^2)}=\frac{t_2}{2(1+t_{2}^2)},\]
 and since $t_1,t_2\in (-1,1)$ this would imply $t_1=t_2$, a contradiction.
 In fact we can show $\phi_{t_1}(\DD)\cap\phi_{t_2}(\DD)=\varnothing$ if and only if $t_1\neq t_2$.
\end{example}
In the infinitesimal case consider the following:
\begin{example}\label{eg:cce_icg_nu}
    For the diamond there exist complex geodesics with $\phi(1)=\psi(1)$, $\phi'(1)=\psi'(1)$ and $\phi(\DD)\neq\psi(\DD)$. 
  \begin{equation}
\begin{aligned}
\phi(\lambda) &= \left(\frac{\lambda}{20}, \frac{19\lambda}{20}\right),\\
\psi(\lambda) &= \left(\frac{1}{5}\left(\lambda-\frac{1}{2}\right)
\left(1-\frac{\lambda}{2}\right),
\mu(\lambda+\nu)(1+\nu\lambda)\right),\\
\mu &= \frac{15+\sqrt{209}}{40}, \qquad
\nu = \frac{15-\sqrt{209}}{4}.
\end{aligned}
\end{equation}
Both $\phi$ and $\psi$ are complex geodesics of the diamond. This follows from Theorem~\ref{thm:geocce}, and one can verify $\phi(1)=\psi(1)$ and $\phi'(1)=\psi'(1)$.

We claim that
$
\phi(\co{\DD})\cap
\psi(\co{\DD})
=\{\phi(1)\}.
$
Indeed, since $\phi_2=19\phi_1$, any point in
$\phi(\co{\DD})\cap
\psi(\co{\DD})$
must satisfy $\psi_2(\lambda)=19\psi_1(\lambda)$
for some $\lambda\in\co{\DD}$.
Using
\[
\mu\nu=\frac{1}{10},
\qquad
\mu(1+\nu^2)=\frac34,
\]
we obtain
\[
\begin{aligned}
\psi_2(\lambda)-19\psi_1(\lambda)
&=\mu(\lambda+\nu)(1+\nu\lambda)\\
&\quad-\frac{19}{5}
\left(\lambda-\frac12\right)
\left(1-\frac{\lambda}{2}\right)\\
&=2(\lambda-1)^2.
\end{aligned}
\]
Hence $\psi_2(\lambda)=19\psi_1(\lambda)$
if and only if $\lambda=1$.
Since $\phi(1)=\psi(1)$, it follows that
\[
\phi(\co{\DD})\cap
\psi(\co{\DD})
=\left\{\left(\frac1{20},\frac{19}{20}\right)\right\}.
\]
In particular,
\[
\phi(\DD)\cap\psi(\DD)=\varnothing.
\]

\end{example}
We do not know whether, for the diamond, complex geodesics passing through an interior point and a boundary point are unique. Nevertheless, we prove the following.
\begin{proposition}\label{prop:cce_ug_inp_bp}
 Let $\cce\subset\CC^n$ be a convex complex ellipsoid with $p_j\in\{1/2,1\}$ for all $j$ and let $\phi,\psi:\DD\to \cce$ be two complex geodesics such that $\phi(0)=\psi(0)$, $\phi(\sigma)=\psi(\sigma)$ and $\phi'(\sigma)=\psi'(\sigma)$ for some $\sigma\in \co{\DD}\setminus\{0\}$. Then $\phi\equiv\psi$.
 \end{proposition}
 If $p_j=1$ for all $j$, then $\cce$ is the unit ball $\BB^n\subset\CC^n$, for which it is already known that any two distinct points of $\co{\BB^n}$ lie on a unique complex geodesic.  
 
 Before going into the proof of the above proposition we need some preparatory results.
\begin{lemma}\label{lem:JWC}
 Let $f\in\OO(\DD,\DD)$ be such that $f(0)=0$ and $f(\sigma)=\sigma$, and suppose $f$ extends holomorphically to a neighborhood of $\sigma\in\ole$ with $f'(\sigma)=1$. Then $f\equiv \id_{\DD}$.
\end{lemma}
\begin{proof}
   This is a consequence of the Julia--Wolff--Carathéodory theorem, see \cite[Corollary~1.2.10]{Aba 1989}.
\end{proof}
\begin{proposition}\label{prop:uni}
Let $D\subset\CC^{n}$ be a convex domain, and let
$\phi,\psi:\DD\to D$ be complex geodesics. Let
$\sigma\in\partial\DD$. Assume that $\phi$ and $\psi$ extend
holomorphically to a neighborhood of $\sigma$ in
$\co{\DD}$, and that
$
\phi(0)=\psi(0),\,
\phi(\sigma)=\psi(\sigma)=b\in\partial D,
$
and
$
\phi'(\sigma)=\psi'(\sigma).
$
 Moreover, assume that $\phi$ admits a holomorphic left inverse 
$F:D\to\DD$, that is,
$
F\circ\phi=\id_{\DD},
$
and that $F$ admits a holomorphic extension to a neighborhood of
$b$.

Then for every $t\in[0,1]$, the map $\eta_t:\DD\to D$
\[
\eta_t(\lambda):=t\phi(\lambda)+(1-t)\psi(\lambda)
\]
is a complex geodesic of $D$.
\end{proposition}
\begin{proof}
Since $D$ is convex and $\phi(\lambda),\psi(\lambda)\in D$ for every
$\lambda\in\DD$, we have
\[
\eta_t(\lambda)=t\phi(\lambda)+(1-t)\psi(\lambda)\in D.
\]

Let
\[
g_t:=F\circ\eta_t:\DD\to\DD.
\]
Then
\[
g_t(0)
=
F(\eta_t(0))
=
F(\phi(0))
=
0,
\]
because $F\circ\phi=\id_{\DD}$.

Moreover, since $\eta_t(\sigma)=b$ and $F$ extends holomorphically, in particular continuously, to a neighborhood of $b$,
we get
\[
g_t(\sigma)
=
F(b)
=
\lim_{\lambda\to\sigma}F(\phi(\lambda))
=
\lim_{\lambda\to\sigma}\lambda
=
\sigma.
\]
Since $F$ extends holomorphically near $b$, we can compute
the boundary derivative:
\[
g_t'(\sigma)
=
dF_b\bigl(\eta_t'(\sigma)\bigr).
\]
But
\[
\eta_t'(\sigma)
=
t\phi'(\sigma)+(1-t)\psi'(\sigma)
=
\phi'(\sigma),
\]
because $\phi'(\sigma)=\psi'(\sigma)$. Therefore
\[
g_t'(\sigma)
=
dF_b\bigl(\phi'(\sigma)\bigr)
=
(F\circ\phi)'(\sigma)
=
1.
\]

Thus, $g_t:\DD\to\DD$ satisfies
\[
g_t(0)=0,\qquad g_t(\sigma)=\sigma,\qquad g_t'(\sigma)=1.
\]
By Lemma~\ref{lem:JWC}, it follows that
\[
g_t(\lambda)=\lambda,\qquad \lambda\in\DD.
\]
Hence
\[
F\circ\eta_t=\id_{\DD}.
\]
Therefore, $\eta_t$ admits a holomorphic left inverse. Consequently,
$\eta_t$ is a complex geodesic of $D$.
\end{proof}

For convex complex ellipsoids with $p_j\in\{1/2,1\}$, the complex geodesics are given explicitly by Theorem~\ref{thm:geocce}, and they extend holomorphically across $\co{\DD}$. The following lemma will be crucial in the proofs below.

\begin{lemma}\label{lem:cce_mid_geo}
Let $\cce\subset\CC^n$ be a convex complex ellipsoid with $p_j\in\{1/2,1\}$ for all $j$, and let $\phi,\psi:\DD\to\cce$ be complex geodesics such that $\eta:=\tfrac12(\phi+\psi)$ is also a complex geodesic of $\cce$. Then for every $j\in\{1,\dots,n\}$ there exist $\theta_j\in\DD$ and a polynomial $R_j$ of degree at most $2$ such that
\begin{equation}\label{eq:mid_geo_form}
\psi_j(\lambda)-\phi_j(\lambda)=\frac{R_j(\lambda)}{(1-\overline{\theta_j}\lambda)^2},
\end{equation}
and $R_j\equiv 0$ whenever $p_j=1$.

In particular, if for every $j$ the function $\psi_j-\phi_j$ has at least three zeros in $\co\DD$, counted with multiplicity, then $\phi\equiv\psi$.
\end{lemma}

\begin{proof}
For $\theta\in\co\DD$ and $r\in\{0,1\}$ put
\[
B_\theta(\lambda):=\frac{\lambda-\theta}{1-\overline{\theta}\lambda},\qquad
P_{\theta,r}(\lambda):=(\lambda-\theta)^{r}(1-\overline{\theta}\lambda)^{2-r}.
\]
If $\theta\in\DD$, then $B_\theta$ is holomorphic near $\co\DD$ and has exactly one zero in $\DD$, a simple zero at $\theta$. If $|\theta|=1$, then $B_\theta\equiv-\theta$ is a unimodular constant.

By Theorem~\ref{thm:geocce} and the remark following it, every component $\phi_j\not\equiv0$ with $p_j=\tfrac12$ can be written as
\[
\phi_j(\lambda)=\frac{a_j\,P_{\alpha_j,r_j}(\lambda)}{(1-\overline{\alpha_0}\lambda)^2},
\]
with $a_j\neq0$, $\alpha_0\in\DD$, $r_j\in\{0,1\}$, $\alpha_j\in\DD$ if $r_j=1$ and $\alpha_j\in\co\DD$ if $r_j=0$. The same holds for $\psi$, with parameters $c_j,\beta_0,\beta_j,s_j$ in place of $a_j,\alpha_0,\alpha_j,r_j$.

Since $p_j\in\{1/2,1\}$, the maps $\phi$, $\psi$ and $\eta$ extend holomorphically to a neighborhood of $\co\DD$, and they map $\ole$ into $\partial\cce$. 

Fix $\xi\in\ole$ and write $u_j=\phi_j(\xi)$ and $v_j=\psi_j(\xi)$. Then
\[
\Bigl|\frac{u_j+v_j}{2}\Bigr|^{2}=\frac{|u_j|^2+|v_j|^2}{2}-\Bigl|\frac{u_j-v_j}{2}\Bigr|^2\quad(p_j=1),\qquad
\Bigl|\frac{u_j+v_j}{2}\Bigr|\le\frac{|u_j|+|v_j|}{2}\quad(p_j=\tfrac12).
\]
Therefore
\[
1=\sum_{j=1}^n|\eta_j(\xi)|^{2p_j}\le\sum_{j=1}^n\frac{|\phi_j(\xi)|^{2p_j}+|\psi_j(\xi)|^{2p_j}}{2}=1,
\]
so equality holds in each of the $n$ inequalities above. Consequently, for every $\xi\in\ole$:
\begin{itemize}
\item if $p_j=1$, then $\phi_j(\xi)=\psi_j(\xi)$;
\item if $p_j=\tfrac12$, then $|u_j+v_j|=|u_j|+|v_j|$. Since $|u_j+v_j|^2=|u_j|^2+|v_j|^2+2\operatorname{Re}(u_j\overline{v_j})$ and $(|u_j|+|v_j|)^2=|u_j|^2+|v_j|^2+2|u_j\overline{v_j}|$, this means $u_j\overline{v_j}\in[0,\infty)$. In particular,
\begin{equation}\label{eq:mid_geo_real}
\phi_j\,\overline{\psi_j}=\psi_j\,\overline{\phi_j}\qquad\text{on }\ole .
\end{equation}
\end{itemize}

\noindent\emph{The case $p_j=1$.} The functions $\phi_j$ and $\psi_j$ are holomorphic near $\co\DD$ and coincide on $\ole$. Hence $\phi_j\equiv\psi_j$, and \eqref{eq:mid_geo_form} holds with $R_j\equiv0$ and $\theta_j=0$.

\noindent\emph{The case $p_j=\tfrac12$, and $\phi_j\equiv0$ or $\psi_j\equiv0$.} If both vanish identically, take $R_j\equiv0$. If $\phi_j\equiv0\not\equiv\psi_j$, write $\psi_j=c_jP_{\beta_j,s_j}/(1-\overline{\beta_0}\lambda)^2$. Then \eqref{eq:mid_geo_form} holds with $\theta_j=\beta_0$ and $R_j=c_jP_{\beta_j,s_j}$. The case $\psi_j\equiv0\not\equiv\phi_j$ is symmetric.

\noindent\emph{The case $p_j=\tfrac12$, and $\phi_j\not\equiv0\not\equiv\psi_j$.} Write
\[
\phi_j=\frac{a_j\,P_{\alpha_j,r_j}}{(1-\overline{\alpha_0}\lambda)^2},\qquad
\psi_j=\frac{c_j\,P_{\beta_j,s_j}}{(1-\overline{\beta_0}\lambda)^2},
\]
with $a_j,c_j\neq0$ and $\alpha_0,\beta_0\in\DD$. For $\xi\in\ole$ and $\theta\in\co\DD$ we have $\overline{1-\overline{\theta}\xi}=\overline{\xi}(\xi-\theta)$ and $\overline{\xi-\theta}=\overline{\xi}(1-\overline{\theta}\xi)$. Hence
\[
\overline{P_{\theta,r}(\xi)}=\overline{\xi}^{\,2}(1-\overline{\theta}\xi)^{r}(\xi-\theta)^{2-r}=\overline{\xi}^{\,2}\,P_{\theta,r}(\xi)\,B_\theta(\xi)^{2(1-r)},
\qquad
\overline{(1-\overline{\alpha_0}\xi)^2}=\overline{\xi}^{\,2}(\xi-\alpha_0)^2 .
\]
Thus, on $\ole$,
\[
\overline{\phi_j}=\frac{\overline{a_j}\,P_{\alpha_j,r_j}B_{\alpha_j}^{2(1-r_j)}}{(\lambda-\alpha_0)^2},\qquad
\overline{\psi_j}=\frac{\overline{c_j}\,P_{\beta_j,s_j}B_{\beta_j}^{2(1-s_j)}}{(\lambda-\beta_0)^2}.
\]
Substitute these into \eqref{eq:mid_geo_real}. The polynomial $P_{\alpha_j,r_j}P_{\beta_j,s_j}\not\equiv0$ has only finitely many zeros, so it can be cancelled. We obtain
\[
a_j\overline{c_j}\,\bigl(B_{\alpha_0}B_{\beta_j}^{\,1-s_j}\bigr)^2=\overline{a_j}c_j\,\bigl(B_{\beta_0}B_{\alpha_j}^{\,1-r_j}\bigr)^2
\]
on $\ole$ minus a finite set. Both sides are rational functions, so the identity holds on all of $\CC$. Both sides are holomorphic on $\DD$, so we may compare their zeros in $\DD$. After halving the multiplicities, the zero divisors in $\DD$ satisfy
\[
[\alpha_0]+S_\beta=[\beta_0]+S_\alpha,\qquad
S_\alpha:=\begin{cases}[\alpha_j], & r_j=0,\ \alpha_j\in\DD,\\ 0, & \text{otherwise},\end{cases}
\]
and similarly for $S_\beta$. Each of $S_\alpha$, $S_\beta$ has degree at most one. So either $\alpha_0=\beta_0$, or $S_\alpha=[\alpha_0]$ and $S_\beta=[\beta_0]$.
\begin{itemize}
\item If $\alpha_0=\beta_0$, then
\[
\psi_j-\phi_j=\frac{c_jP_{\beta_j,s_j}-a_jP_{\alpha_j,r_j}}{(1-\overline{\alpha_0}\lambda)^2},
\]
which is \eqref{eq:mid_geo_form} with $\theta_j=\alpha_0$.
\item Otherwise $r_j=s_j=0$, $\alpha_j=\alpha_0$ and $\beta_j=\beta_0$. Then $\phi_j\equiv a_j$ and $\psi_j\equiv c_j$ are constant, and \eqref{eq:mid_geo_form} holds with $\theta_j=0$ and $R_j\equiv c_j-a_j$.
\end{itemize}

\smallskip
\noindent\textbf{The final assertion.} Since $\theta_j\in\DD$, the factor $(1-\overline{\theta_j}\lambda)^{-2}$ is holomorphic and zero-free on the disc $\{|\lambda|<1/|\theta_j|\}\supset\co\DD$. Hence $\psi_j-\phi_j$ and $R_j$ have the same zeros, with the same multiplicities, on $\co\DD$. A polynomial of degree at most $2$ with at least three zeros counted with multiplicity vanishes identically. So $R_j\equiv0$ for every $j$, that is, $\phi\equiv\psi$.
\end{proof}

 Now we give the proof of Proposition~\ref{prop:cce_ug_inp_bp}.
 \begin{proof}[Proof of Proposition~\ref{prop:cce_ug_inp_bp}]
If $\sigma\in\DD\setminus\{0\}$, it is shown in \cite{Jar-Pfl-Zei 1993} that $\phi(0)=\psi(0)$ and $\phi(\sigma)=\psi(\sigma)$ already force $\phi\equiv\psi$.
 It therefore remains to treat the case $\sigma\in\ole$ with $\psi(\sigma)=\phi(\sigma)$ and $\phi'(\sigma)=\psi'(\sigma)$.

Define
\[
\eta(\lambda):=\frac{\phi(\lambda)+\psi(\lambda)}{2},\qquad \lambda\in\DD.
\]
Since $\cce$ is convex, $\eta(\DD)\subset\cce$. Since $p_j\in\{1/2,1\}$, the geodesics $\phi$ and $\psi$ extend holomorphically to a neighborhood of $\co\DD$, and $\phi(\sigma)=\psi(\sigma)\in\partial\cce$. By Theorem~\ref{thm:extn_li_cce_1_2} (below), $\phi$ has a left inverse that extends holomorphically through $\phi(\ole)$. Hence, by Proposition~\ref{prop:uni}, $\eta$ is a complex geodesic.

For every $j$, the function $\psi_j-\phi_j$ vanishes at $0$ and vanishes to order at least $2$ at $\sigma\neq0$, so it has at least three zeros in $\co\DD$, counted with multiplicity. By Lemma~\ref{lem:cce_mid_geo}, $\phi\equiv\psi$.
\end{proof}

\subsection{Existence of a left inverse and application of Zwonek's lemma}
In this subsection we show the existence of holomorphically extendable left inverses of complex geodesics of $\cce$ with $p_j\in\{1/2,1\}$. Using such a left inverse and Zwonek's lemma (Lemma~\ref{lem:zwobk}), we show that, under a Burns--Krantz map of $\cce$ at $b\in \partial \cce$, a complex geodesic through $b$ is mapped to a complex geodesic. 

 To begin with, recall that \cite[Proposition~3.1]{Zwo 2025b} uses the dual map of a complex geodesic in a bounded convex domain to extend its Lempert left inverse holomorphically. 
 \begin{proposition}[{\cite[Proposition~3.1]{Zwo 2025b}}]
Let $D$ be a bounded convex domain in $\mathbb C^n$ and let $\phi$ be a complex geodesic in $D$, $\widetilde{\phi}$ its dual and $F$ the Lempert left inverse to $\phi$ determined by $\widetilde{\phi}$ (which in particular means that it satisfies the equation (\ref{eq:dual_lem_inv})). Assume further that $$\langle \phi', \widetilde{\phi}\rangle\equiv C>0$$ on $\mathbb D$. Fix $z_0\in\partial D$ and let $\lambda_0\in\overline{\mathbb D}$ be such that $\lambda_0$ belongs to a cluster set $F^{c}(z_0)$. Assume that $\phi$ and $\widetilde{\phi}$ extend holomorphically through $\lambda_0$ (which is trivially satisfied if $\lambda_0\in\mathbb D$). Then $F$ extends holomorphically through $z_0$.
\end{proposition}
Recall the cluster set of $F:D\to \DD$ at $z_{0}\in\partial{D}$ is: $$F^{c}(z_0)=\left \{\lim_{k\to \infty} F(z_k): D\ni \{z_k\}\to z_0\text{ and the limit exists}\right \}$$
 We remark that the above proposition is not correct as stated. Consider the following example for the diamond. The cluster set of the Lempert left inverse $F(z_1,z_2)=z_1/(1-z_2)$ of the complex geodesic $\phi(\lambda)=(\lambda/2,1/2)$ (which is determined by the dual map of $\phi$ from Equation~(\ref{eq:dual_lem_inv})), at the point $z_0=(0,1)\in \partial \triangle$ is $F^{c}(z_0)=\co{\DD}$. Indeed, let $w\in {\DD}$, take a real sequence $t_k\searrow 0$, so that $z_{k}=(wt_k,1-t_k)\in {\triangle}$, $z_k\to z_0$ and $F(z_k)=w$. Hence $\DD\subset F^{c}(z_0)$. Since cluster sets are closed and $F(\triangle)\subset \DD$, we get $F^{c}(z_0)=\co{\DD}$. So \cite[Proposition~3.1]{Zwo 2025b} would imply that $F$ extends holomorphically near $(0,1)$ which is not possible as $F$ is not even continuous at $z_0$.  

Nevertheless, the idea behind the proof of that proposition is enough for our purpose to obtain the following.

\begin{theorem}[{\cite{Zwo 2025b}}]\label{thm:zwo_dual}
    Let $D$ be a bounded convex domain in $\CC^n$, let $\phi:\DD\to D$ be a complex geodesic, let $\widetilde{\phi}$ be a dual map of $\phi$, and let $F$ be the Lempert left inverse of $\phi$ determined by $\widetilde{\phi}$ (in particular, $F$ satisfies \eqref{eq:dual_lem_inv}). Assume further that 
    \begin{equation}\label{eq:pos_crit_dual}
        \langle \phi'(\lambda),\widetilde{\phi}(\lambda)\rangle \equiv C>0\text{ on } \DD.
    \end{equation}  
    Let $\lambda_0\in \ole$ be such that $\phi$ and $\widetilde{\phi}$ extend holomorphically through $\lambda_0$. Then $F$ extends holomorphically through $z_0:=\phi(\lambda_0)\in \partial D$.
\end{theorem}
\begin{proof}
    Let $U$ be an open neighborhood of $\lambda_0$ to which $\phi$ and $\widetilde{\phi}$ extend holomorphically, and define the holomorphic function
    \[
        \Phi(z,\lambda):=\langle z-\phi(\lambda),\widetilde{\phi}(\lambda)\rangle,\qquad (z,\lambda)\in\CC^n\times U.
    \]
    By Equation~\eqref{eq:pos_crit_dual} and continuity, $\langle\phi'(\lambda_0),\widetilde{\phi}(\lambda_0)\rangle=C$. Hence $\Phi(z_0,\lambda_0)=0$ and
    \[
        \frac{\partial\Phi}{\partial\lambda}(z_0,\lambda_0)=-\langle\phi'(\lambda_0),\widetilde{\phi}(\lambda_0)\rangle+\langle z_0-\phi(\lambda_0),\widetilde{\phi}'(\lambda_0)\rangle=-C\neq0.
    \]
    Now by application of the holomorphic implicit function theorem, $F$ extends holomorphically through $z_0$. 
    
    Indeed, the implicit function theorem gives neighborhoods $V\ni z_0$, $W\ni\lambda_0$ and a holomorphic map $\widehat F:V\to W$ such that, for $z\in V$, $\widehat F(z)$ is the unique zero of $\Phi(z,\cdot)$ in $W$. Shrinking $V$, we may assume that $V$ is a ball, so $V\cap D$ is connected.
For $\mu\in\partial\DD\cap W$, the vector $\mu\overline{\widetilde{\phi}(\mu)}\neq0$ is an outward normal to $D$ at $\phi(\mu)$ \cite[Section 2.6.3]{Aba 1989}, hence $\operatorname{Re}\bigl(\overline{\mu}\,\Phi(z,\mu)\bigr)<0$ for $z\in D$. In particular $\widehat F(V\cap D)\cap\partial\DD=\varnothing$. Since $\widehat F(\phi(\lambda))=\lambda\in\DD$ for $\lambda\in\DD$ close to $\lambda_0$, connectedness gives $\widehat F(V\cap D)\subset\DD$. As $F(z)$ is the unique zero of $\Phi(z,\cdot)$ in $\DD$, we get $\widehat F=F$ on $V\cap D$, so $F$ extends holomorphically through $z_0$.

\end{proof}
Using Theorem~\ref{thm:zwo_dual}, we prove that complex geodesics of $\cce$ with $p_j\in\{1/2,1\}$ admit left inverses that extend holomorphically. 

Recall from Theorem~\ref{thm:geocce} that every complex geodesic $\psi:\DD\to\cce$ with $p_j\geq 1/2$ and $\psi_{j}\not \equiv 0$ has the form
\begin{equation}
    \psi_j(\lambda)=a_j\left(\frac{\lambda-\alpha_j}{1-\overline{\alpha_j}\lambda}\right)^{r_j}\left(\frac{1-\overline{\alpha_j}\lambda}{1-\overline{\alpha_0}\lambda}\right)^{1/p_j},\quad j=1,\dots,n,
\end{equation}
where $r_j\in\{0,1\}$ and the constants satisfy the relations in Theorem~\ref{thm:geocce}. When $p_j\in\{1/2,1\}$, $\psi$ extends holomorphically across $\co{\DD}$.

The dual maps of these complex geodesics were described explicitly in \cite{Jar-Pfl-Zei 1993}: a dual map of $\psi$ is $\widetilde{\psi}:\DD\to\CC^n$, given by
\begin{equation}\label{eq:dual_half_one}
    \widetilde{\psi}_j(\lambda)=p_j|a_j|^{2p_j}\,\frac{(\lambda-\alpha_j)(1-\overline{\alpha_j}\lambda)}{\psi_j(\lambda)},\quad j=1,\dots,n.
\end{equation}
If $\psi_j\equiv 0$ for some $j$, we set $\widetilde{\psi}_j\equiv 0$.
For $p_j\in\{1/2,1\}$, the singularities of $\widetilde{\psi}$ in $\co{\DD}$ are removable, so $\widetilde{\psi}$ also extends holomorphically to a neighborhood of $\co{\DD}$.

We begin by showing that the positivity condition \eqref{eq:pos_crit_dual} is satisfied on every convex complex ellipsoid $\cce$.

\begin{lemma}\label{lem:pos_dual_geo}
Let $\cce\subset\CC^n$ be a convex complex
ellipsoid. Let
$\psi:\DD\to\cce$ be a complex geodesic as in Theorem~\ref{thm:geocce}.
If $\widetilde{\psi}$ is the dual map of $\psi$ given by \eqref{eq:dual_half_one}, then
\[
\langle\psi'(\lambda),\widetilde{\psi}(\lambda)\rangle
\equiv C>0,\qquad \lambda\in\DD.
\]
\end{lemma}

\begin{proof}
By Theorem~\ref{thm:geocce}, the components of
$\psi$ are given by
\[
\psi_j(\lambda)
=
a_j
\left(
\frac{\lambda-\alpha_j}
{1-\overline{\alpha_j}\lambda}
\right)^{r_j}
\left(
\frac{1-\overline{\alpha_j}\lambda}
{1-\overline{\alpha_0}\lambda}
\right)^{1/p_j},
\]
where $a_j\neq0$, $\alpha_0\in\DD$,
$\alpha_j\in\co{\DD}$, and
$\alpha_j\in\DD$ whenever $r_j=1$.

Logarithmic differentiation yields
\begin{align*}
\frac{\psi_j'(\lambda)}{\psi_j(\lambda)}
&=
r_j\frac{1-|\alpha_j|^2}
{(\lambda-\alpha_j)
(1-\overline{\alpha_j}\lambda)}\\
&\quad+
\frac{1}{p_j}
\left(
\frac{\overline{\alpha_0}}
{1-\overline{\alpha_0}\lambda}
-
\frac{\overline{\alpha_j}}
{1-\overline{\alpha_j}\lambda}
\right).
\end{align*}

Multiplying by $\widetilde{\psi}_j(\lambda)\psi_j(\lambda)$
and simplifying, we obtain
\[
\psi_j'(\lambda)\widetilde{\psi}_j(\lambda)
=
p_j|a_j|^{2p_j}r_j(1-|\alpha_j|^2)
+
|a_j|^{2p_j}
\frac{(\lambda-\alpha_j)
(\overline{\alpha_0}-\overline{\alpha_j})}
{1-\overline{\alpha_0}\lambda}.
\]
The identity extends across any removable
singularities.

Put
\[
\omega_j:=|a_j|^{2p_j}>0,
\qquad
S:=\sum_{j=1}^n \omega_j.
\]
Recall the parameter relations from
Theorem~\ref{thm:geocce}:
\[
\alpha_0=\sum_{j=1}^n \omega_j\alpha_j,
\qquad
1+|\alpha_0|^2
=
\sum_{j=1}^n \omega_j(1+|\alpha_j|^2).
\]
In particular,
\[
\sum_{j=1}^n \omega_j|\alpha_j|^2
=
1+|\alpha_0|^2-S.
\]

Consequently,
\begin{align*}
&\sum_{j=1}^n
\omega_j(\lambda-\alpha_j)
(\overline{\alpha_0}-\overline{\alpha_j})\\
&\quad=
\lambda\left(
S\overline{\alpha_0}
-\sum_{j=1}^n \omega_j\overline{\alpha_j}
\right)
-\overline{\alpha_0}
\sum_{j=1}^n \omega_j\alpha_j
+\sum_{j=1}^n \omega_j|\alpha_j|^2\\
&\quad=
\lambda(S-1)\overline{\alpha_0}
-|\alpha_0|^2
+1+|\alpha_0|^2-S \\
&\quad=
(1-S)(1-\overline{\alpha_0}\lambda).
\end{align*}

It follows that
\begin{align*}
\langle\psi'(\lambda),\widetilde{\psi}(\lambda)\rangle
&=
\sum_{j=1}^n\psi_j'(\lambda)\widetilde{\psi}_j(\lambda)\\
&=
1-S+
\sum_{j=1}^n
p_j\omega_jr_j(1-|\alpha_j|^2).
\end{align*}
Thus $\langle\psi',\widetilde{\psi}\rangle$ is a real constant.

To prove its positivity, observe that the
parameter relations imply
\begin{align*}
\sum_{j=1}^n \omega_j|\alpha_j-\alpha_0|^2
&=
\sum_{j=1}^n \omega_j|\alpha_j|^2
-2|\alpha_0|^2+S|\alpha_0|^2\\
&=
(1-S)(1-|\alpha_0|^2).
\end{align*}
Since $\alpha_0\in\DD$, we obtain
\[
1-S
=
\frac{\displaystyle\sum_{j=1}^n
\omega_j|\alpha_j-\alpha_0|^2}
{1-|\alpha_0|^2}
\geq0.
\]

Moreover, each term
\[
p_j\omega_jr_j(1-|\alpha_j|^2)
\]
is nonnegative, and is strictly positive
whenever $r_j=1$.

Therefore, if $r_j=1$ for at least one $j$,
we immediately obtain
\[
\langle\psi',\widetilde{\psi}\rangle>0.
\]

It remains to consider the case when
$r_j=0$ for every $j$. In this case,
\[
\langle\psi',\widetilde{\psi}\rangle=1-S \geq0.
\]
Suppose that $S=1$. Then
\[
\sum_{j=1}^n \omega_j|\alpha_j-\alpha_0|^2=0.
\]
Since $\omega_j>0$ for every $j$, it follows that
\[
\alpha_j=\alpha_0
\qquad\text{for all }j.
\]
Consequently,
\[
\psi_j(\lambda)
=
a_j
\left(
\frac{1-\overline{\alpha_0}\lambda}
{1-\overline{\alpha_0}\lambda}
\right)^{1/p_j}
=a_j.
\]
Thus $\psi$ is constant, contradicting the
assumption that $\psi$ is a complex geodesic.
Hence $S<1$, which again yields
$\langle\psi',\widetilde{\psi}\rangle>0$.

If $\psi_j\equiv 0$ for some $j$, then $\widetilde{\psi}_j\equiv 0$, so the $j$-th term of $\langle\psi',\widetilde{\psi}\rangle$ vanishes. Applying the same computation to the indices $k$ with $\psi_k\not\equiv 0$ shows that the result still holds.

We conclude that
\[
\langle\psi'(\lambda),\widetilde{\psi}(\lambda)\rangle
=
1-S+
\sum_{j=1}^n
p_j\omega_jr_j(1-|\alpha_j|^2)
=:C>0
\]
for every $\lambda\in\DD$.
\end{proof}
Now we apply Theorem~\ref{thm:zwo_dual} and show that complex geodesics of $\cce$ with $p_j\in\{1/2,1\}$ have a Lempert left inverse that extends holomorphically.
\begin{theorem}\label{thm:extn_li_cce_1_2}
     Let $\cce\subset\CC^n$ be a convex complex ellipsoid with $p_{j}\in\{1/2,1\}$ for all $j\in \{1,\dots,n\}$, and let $\psi:\DD\to \cce$ be a complex geodesic. Then there exists a Lempert left inverse of $\psi$ which extends holomorphically through $\psi(\ole)\subset \partial \cce$.
\end{theorem}
\begin{proof}
    Observe that $\cce$ is a bounded convex domain. The complex geodesic $\psi$ and the dual map $\widetilde{\psi}$ given by \eqref{eq:dual_half_one} both extend holomorphically to a neighborhood of $\co{\DD}$. Let $F$ be the Lempert left inverse determined by $\widetilde{\psi}$ as in \eqref{eq:dual_lem_inv}. By Lemma~\ref{lem:pos_dual_geo} we have
    \[\langle \psi'(\lambda),\widetilde{\psi}(\lambda)\rangle \equiv C>0\text{ on } \DD.\]
    Fix $\zeta_0\in \ole$. Since $\psi$ and $\widetilde{\psi}$ extend holomorphically through $\zeta_0$, all the hypotheses of Theorem~\ref{thm:zwo_dual} are satisfied with $\lambda_0=\zeta_0$. Hence $F$ extends holomorphically through $\psi(\zeta_0)$. As $\zeta_0\in \ole$ was arbitrary, $F$ extends holomorphically through $\psi(\ole)\subset \partial \cce$.
\end{proof}
We use Zwonek's lemma (Lemma~\ref{lem:zwobk}) to show that under a Burns--Krantz map, a complex geodesic is mapped to a complex geodesic.
\begin{proposition}\label{prop:psi_cg}
    Let $\cce\subset \CC^n$ be a convex complex ellipsoid with $p_j\in \{1/2,1\}$ for all $j\in \{1,\dots,n\}$. Let $\phi:\DD\to \cce$ be a complex geodesic with $\phi(1)=b$ and $f:\cce\to \cce$ be a Burns--Krantz map at $b\in \partial \cce$. Then $\psi:=f\circ\phi:\DD\to \cce$ is a complex geodesic with the same Lempert left inverse $G:\cce\to \DD$ as $\phi$, where $G$ extends holomorphically to a neighborhood of $\phi(\ole)\subset \partial \cce$. 
\end{proposition}
\begin{proof}
    As complex geodesics of $\cce$ extend holomorphically to a neighborhood of $\co{\DD}$, in particular $\phi$ is Lipschitz near $1$. By Theorem~\ref{thm:extn_li_cce_1_2} there exists a Lempert left inverse $G:\cce\to \DD$ of $\phi$ that extends holomorphically to a neighborhood of $\phi(\ole)\subset\partial\cce$. In particular $G$ is Lipschitz at $b\in \partial\cce$. 
    
    Hence, by Zwonek's lemma (Lemma~\ref{lem:zwobk}), we get that $\psi=f\circ \phi$ is a complex geodesic with the same Lempert left inverse $G$.  
\end{proof}
\subsection{Equality of derivatives at the boundary}
In this subsection we show that $\psi=f\circ\phi$ and $\phi$ agree up to the third derivative at the boundary point $1\in \ole$. This will imply that $\psi\equiv \phi$.

We begin by showing that if two complex geodesics of a convex domain share a common Lempert left inverse then their convex combination is also a complex geodesic. Here the fact that the left inverse is Lempert plays a crucial role. For example, $\widehat{F}(z_1,z_2)=4z_1z_2$ is a common non-Lempert left inverse for complex geodesics of the diamond $\phi(\lambda)=(\lambda/2,1/2)$ and $\psi(\lambda)=(-\lambda /2,-1/2)$, but the map $(\phi+\psi)/2\equiv 0$ is not a complex geodesic.
\begin{lemma}\label{lem:cmn_lem_inv}
    Let $D\subset\CC^n$ be a bounded convex domain. Let $\phi,\psi:\DD\to D$ be two complex geodesics with a common Lempert left inverse $G:D\to \DD$. Then $\eta_{t}(\lambda):=t\phi(\lambda)+(1-t)\psi(\lambda)$ is a complex geodesic with the same Lempert left inverse $G$ for any $t\in [0,1]$.
\end{lemma}
\begin{proof}
    Since $D$ is convex, $\eta_t(\DD)\subset D$ for $t\in [0,1]$.
    
    As $G$ is a Lempert left inverse (see Equation~\eqref{eq:lem_inv_hp}), for $\lambda\in\DD$
    $$G^{-1}(\lambda)\cap D = D\cap H_{\lambda}$$
    where $H_{\lambda}$ is an affine complex hyperplane. Furthermore, $G\circ\phi=\id_{\DD}$ and $G\circ\psi=\id_{\DD}$. As a result $\phi(\lambda),\psi(\lambda)\in D\cap H_{\lambda}$. Since $D$ and $H_{\lambda}$ are convex we get $$t\phi(\lambda)+(1-t)\psi(\lambda)\in D\cap H_{\lambda} \text{ for } t\in [0,1].$$
    Consequently, $G(\eta_{t}(\lambda))=\lambda$ for all $\lambda\in \DD$ and any $t\in [0,1]$. In other words for any $t\in [0,1]$, $\eta_{t}$ is a complex geodesic with the same Lempert left inverse $G:D\to \DD$.
\end{proof}
We now show the main result of this subsection. 
\begin{theorem}\label{thm:mod_zwo}
 Let $\cce\subset \CC^n$ be a convex complex ellipsoid with $p_{j}\in\{1/2,1\}$ for all $j\in \{1,\dots,n\}$, let $b \in \partial \cce$, and let
$f: \cce \to \cce$ be a Burns--Krantz map at $b$. Assume that $\phi: \DD \to \cce$ is a complex geodesic such that $\phi(1) = b$, and let $\psi=f\circ\phi$. Then $\phi\equiv\psi$.
\end{theorem}
 
 \begin{proof}
 By Proposition~\ref{prop:psi_cg}, $\psi=f\circ\phi$ is a complex geodesic of $\cce$. Hence, since $p_j\in\{1/2,1\}$, both $\phi$ and $\psi$ extend holomorphically to a neighborhood of $\co\DD$, and so does
\[
h:=\psi-\phi .
\]
Observe that $h(\lambda)=f(\phi(\lambda))-\phi(\lambda)=o(\|\phi(\lambda)-b\|^3)=o(|\lambda-1|^3)$ as $\lambda\to1$, $\lambda\in\DD$.  Since $h$ is holomorphic near $1$,
\begin{equation}\label{eq:order}
h(1)=h'(1)=h''(1)=h'''(1)=0.
\end{equation}
In other words, every component $h_j=\psi_j-\phi_j$ has a zero of order at least $4$ at $\lambda=1$. 

Since by Proposition~\ref{prop:psi_cg}, $\phi$ and $\psi$ have a common Lempert left inverse, Lemma~\ref{lem:cmn_lem_inv} shows that $\tfrac12(\phi+\psi)$ is a complex geodesic and Lemma~\ref{lem:cce_mid_geo} gives $\phi\equiv\psi$.
\end{proof}

\subsection{Proof of the main theorem}
We give the proof of Theorem~\ref{thm:bkcce} by collecting all the above results. 
\begin{proof}
 Let $f:\cce\to\cce$ be a Burns--Krantz map at $b\in\partial\cce$, and let $q\in \cce$ be an arbitrary point.

 Since $\cce$ is $m$-convex, by Mercer's Theorem~\ref{thm:mer_93} there exists a complex geodesic $\phi_{q}:\DD\to \cce$ such that $\phi_{q}(1)=b$ and $\phi_{q}(0)=q$. Since every complex geodesic of $\cce$ with $p_j\in \{1/2,1\}$ extends holomorphically across $\co{\DD}$, so does $\phi_{q}$.

Set $\psi_{q}=f\circ\phi_{q}$. By Proposition~\ref{prop:psi_cg}, $\psi_q$ is also a complex geodesic and by Theorem~\ref{thm:mod_zwo} $\phi_{q}\equiv \psi_{q}$. In other words $f\circ\phi_{q}\equiv \phi_{q}$. In particular, $f(q)=f(\phi_{q}(0))=\phi_{q}(0)=q$. Since $q\in\cce$ was arbitrary, $f\equiv \id_{\cce}$.
\end{proof}
\include{Appendix}
\section*{Declaration on the use of AI tools}
The main result and its proof were obtained without the use of AI tools. The proof of the $m$-convexity of $\cce$ (Proposition~\ref{prop:cce_m_cnv}, Appendix~\ref{Appendix}) was obtained with the help of ChatGPT 5.4 (OpenAI), and Claude (Anthropic) was used to proofread the manuscript for grammar, notation and mathematical errors. AI tools were also used to improve the clarity of the exposition. The author has checked all content and takes full responsibility for the article.

\end{document}

%% file: Appendix.tex
\section{Appendix}\label{Appendix}
We first prove the following estimates, which will be crucial in the proof of Proposition~\ref{prop:cce_m_cnv}.

\begin{lemma}[Circular mean estimate]
\label{lem:circular_mean}
Let $\alpha\geq 1$. Then there exists a constant
$c_\alpha>0$ such that
\begin{equation}\label{eq:circular_mean}
\frac{1}{2\pi}\int_0^{2\pi}
\left(
|a+e^{it}b|^\alpha-|a|^\alpha
\right)\,dt
\geq
c_\alpha |b|^{\max\{2,\alpha\}}
\end{equation}
for all $a,b\in\CC$ satisfying
$|a|\leq 1$ and $|b|\leq 1$.
\end{lemma}

\begin{proof}
The assertion is immediate if $b=0$, so assume
$b\neq 0$. We distinguish two cases.

\smallskip
\noindent
\textbf{Case 1: $|a|\leq 2|b|$.}

By a change of variables,
\[
\frac{1}{2\pi}\int_0^{2\pi}
\left(
|a+e^{it}b|^\alpha-|a|^\alpha
\right)\,dt
=
|b|^\alpha F_\alpha(a/b),
\]
where
\[
F_\alpha(w):=
\frac{1}{2\pi}\int_0^{2\pi}
|w+e^{it}|^\alpha\,dt-|w|^\alpha.
\]

Since $|\zeta|^\alpha$ is strictly subharmonic
for $\alpha\geq 1$, the strict mean-value
inequality implies
\[
F_\alpha(w)>0,\qquad w\in\CC.
\]
Moreover, $F_\alpha$ is continuous. Hence,
by compactness, there exists $c_1>0$, depending only on $\alpha$, such that
\[
F_\alpha(w)\geq c_1
\qquad\text{whenever } |w|\leq 2.
\]
Consequently,
\[
\frac{1}{2\pi}\int_0^{2\pi}
\left(
|a+e^{it}b|^\alpha-|a|^\alpha
\right)\,dt
\geq c_1|b|^\alpha.
\]
Since $|b|\leq 1$ and
$\max\{2,\alpha\}\geq\alpha$, we obtain
\[
\frac{1}{2\pi}\int_0^{2\pi}
\left(
|a+e^{it}b|^\alpha-|a|^\alpha
\right)\,dt
\geq c_1|b|^{\max\{2,\alpha\}}.
\]

\smallskip
\noindent
\textbf{Case 2: $|a|>2|b|$.}

In this case, the closed disc
$\{\zeta\in\CC:|\zeta-a|\leq|b|\}$, centered at $a$ with radius $|b|$, does not
contain the origin. Recall that
\[
\Delta|\zeta|^\alpha
=
\alpha^2|\zeta|^{\alpha-2},
\qquad \zeta\neq 0.
\]
Furthermore, for
$|\zeta-a|\leq |b|$, we have
\[
\frac{|a|}{2}\leq|\zeta|
\leq\frac{3|a|}{2}.
\]
Thus, for some constant $c_2>0$ depending
only on $\alpha$,
\[
\Delta|\zeta|^\alpha
\geq c_2|a|^{\alpha-2}
\]
throughout this disc.

By the Poisson--Jensen formula
(see \cite[Theorem 4.5.1]{Ran 1995}) for the disc of radius $|b|$ centered at $a$ we get
\[
\frac{1}{2\pi}\int_0^{2\pi}
\left(
|a+e^{it}b|^\alpha-|a|^\alpha
\right)\,dt
=
\frac{1}{2\pi}\int_{|\zeta-a|<|b|}
\log\frac{|b|}{|\zeta-a|}\,\Delta|\zeta|^\alpha\,dA(\zeta),
\]
and since $\int_{|\zeta-a|<|b|}\log\frac{|b|}{|\zeta-a|}\,dA(\zeta)=\frac{\pi}{2}|b|^2$, we get
\[
\frac{1}{2\pi}\int_0^{2\pi}
\left(
|a+e^{it}b|^\alpha-|a|^\alpha
\right)\,dt
\geq
\frac{c_2}{4}|b|^2|a|^{\alpha-2}.
\]

If $\alpha\geq 2$, then $|a|>2|b|$
implies
\[
|b|^2|a|^{\alpha-2}
\geq 2^{\alpha-2}|b|^\alpha.
\]

On the other hand, if $1\leq\alpha<2$,
then $|a|\leq 1$, and hence
\[
|a|^{\alpha-2}\geq 1.
\]
Therefore,
\[
|b|^2|a|^{\alpha-2}\geq |b|^2.
\]

Combining the two cases, we obtain, with $c_\alpha:=\min\{c_1,c_2/4\}$,
\[
\frac{1}{2\pi}\int_0^{2\pi}
\left(
|a+e^{it}b|^\alpha-|a|^\alpha
\right)\,dt
\geq
c_\alpha |b|^{\max\{2,\alpha\}},
\]
as required.
\end{proof}

\begin{lemma}[Directional radius estimate]
\label{lem:directional_radius}
Let
\[
\cce=
\left\{
z\in\CC^n:
\rho(z):=\sum_{j=1}^n|z_j|^{2p_j}<1
\right\},
\qquad p_j\geq\frac12,
\]
and put
\[
m:=\max\{2,2p_1,\dots,2p_n\}.
\]
Then there exists a constant $C>0$ such that
\begin{equation}\label{eq:radius_rho}
\delta_{\cce}(z;v)^m
\leq C(1-\rho(z))
\end{equation}
for every $z\in\cce$ and every unit vector
$v\in\CC^n$.
\end{lemma}

\begin{proof}
Put $\kappa_j:=2p_j$ and
$m_j:=\max\{2,\kappa_j\}$. Then
$\kappa_j\geq 1$ and $m_j\leq m$.

Fix $z\in\cce$ and $v\in\CC^n$ with
$\|v\|=1$. Let $r>0$ be such that
\[
z+\lambda v\in\cce
\qquad\text{for all }|\lambda|<r.
\]

For $0<s<r$ we have
\[
\rho(z+se^{it}v)<1,
\qquad t\in[0,2\pi].
\]
Averaging over $t$, we obtain
\[
\frac{1}{2\pi}\int_0^{2\pi}
\rho(z+se^{it}v)\,dt\leq 1.
\]
Consequently,
\begin{align*}
1-\rho(z)
&\geq
\frac{1}{2\pi}\int_0^{2\pi}
\bigl(\rho(z+se^{it}v)-\rho(z)\bigr)\,dt\\
&=
\sum_{j=1}^n
\frac{1}{2\pi}\int_0^{2\pi}
\left(
|z_j+se^{it}v_j|^{\kappa_j}
-|z_j|^{\kappa_j}
\right)\,dt.
\end{align*}

Since the complex disc
$z+\{\lambda v:|\lambda|<r\}$
is contained in $\cce$, we have
$|z_j|\leq 1$ and $s|v_j|\leq 1$
for every $j$.

Therefore, applying
Lemma~\ref{lem:circular_mean} coordinatewise,
we get
\[
1-\rho(z)
\geq
\sum_{j=1}^n c_j|sv_j|^{m_j},
\]
where $c_j:=c_{\kappa_j}>0$ is the constant of Lemma~\ref{lem:circular_mean} for $\alpha=\kappa_j$.

Since $m_j\leq m$ and $s|v_j|\leq 1$,
\[
|sv_j|^{m_j}\geq |sv_j|^m.
\]
Hence, with $c:=\min_j c_j>0$,
\[
1-\rho(z)
\geq
c s^m\sum_{j=1}^n|v_j|^m.
\]

Finally, since $\|v\|=1$ and $m\geq 2$, by Jensen's inequality
applied to the convex function $t\mapsto t^{m/2}$,
we obtain
\[
\sum_{j=1}^{n}|v_j|^m
\geq n^{1-m/2}
\left(\sum_{j=1}^{n}|v_j|^2\right)^{m/2}
= n^{1-m/2}.
\]
Consequently,
\[
1-\rho(z)\geq c\,n^{1-m/2}s^m.
\]

Thus, with $C:=n^{m/2-1}/c$,
\[
s^m\leq C(1-\rho(z)).
\]

Letting $s\uparrow r$ and then taking
the supremum over all admissible radii $r$,
we obtain
\[
\delta_{\cce}(z;v)^m
\leq C(1-\rho(z)).
\]
This completes the proof.
\end{proof}

We now complete the proof of Proposition~\ref{prop:cce_m_cnv}.
\begin{proof}
By Lemma~\ref{lem:directional_radius},
there exists $C>0$ such that
\[
\delta_{\cce}(z;v)^m
\leq C(1-\rho(z))
\]
for all $z\in\cce$ and all unit vectors
$v\in\CC^n$.

Since $2p_j\geq 1$, the defining function
\[
\rho(z)=\sum_{j=1}^n|z_j|^{2p_j}
\]
is Lipschitz continuous on
$\co{\cce}$.
Thus there exists $L>0$ such that
\[
|\rho(z)-\rho(w)|
\leq L\|z-w\|,
\qquad z,w\in\co{\cce}.
\]

For $z\in\cce$, choose $w\in\partial\cce$
such that
\[
\|z-w\|=\delta_{\cce}(z).
\]
Since $\rho(w)=1$, we obtain
\[
1-\rho(z)
=\rho(w)-\rho(z)
\leq L\delta_{\cce}(z).
\]

Combining this with the preceding estimate,
we conclude that
\[
\delta_{\cce}(z;v)^m
\leq CL\delta_{\cce}(z),
\]
and therefore, with $C':=(CL)^{1/m}$,
\[\delta_{\cce}(z;v)
\leq C'\delta_{\cce}(z)^{1/m}.\]
Hence $\cce$ is $m$-convex in the sense
of Definition~\ref{def:m_cvx}.
\end{proof}

%% file: Burns--Krantz_Rigidity_arXiv_.bbl
\begin{thebibliography}{...}
\bibitem[Aba89]{Aba 1989} M. Abate. \emph{Iteration Theory of Holomorphic Maps on Taut Manifolds}. Research and Lecture Notes in Mathematics. Complex Analysis and Geometry. Mediterranean Press, Rende, 1989. MR1098711.
\bibitem[BK94]{Bur-Kra 1994} D. Burns and S.~G. Krantz. Rigidity of holomorphic mappings and a new Schwarz lemma at the boundary. \emph{J. Amer. Math. Soc.}, \textbf{7}(3):661--676, 1994. MR1242454.
\bibitem[Car31]{Cartan1931} H. Cartan. Les fonctions de deux variables complexes et le probl\`eme de la repr\'esentation analytique. \emph{J. Math. Pures Appl. (9)}, \textbf{10}:1--114, 1931.
\bibitem[CHL88]{Cha-Hu-Lee 1988} C.-H. Chang, M.~C. Hu, and H.-P. Lee. Extremal analytic discs with prescribed boundary data. \emph{Trans. Amer. Math. Soc.}, \textbf{310}(1):355--369, 1988. MR0930081.
\bibitem[CS19]{Chr-Sho 2019} A. Christodoulou and I. Short. A hyperbolic-distance inequality for holomorphic maps. \emph{Ann. Acad. Sci. Fenn. Math.}, \textbf{44}(1):293--300, 2019. MR3919138.
\bibitem[FR21]{For-Ron 2021} J.~E. Forn\ae ss and F. Rong. The boundary rigidity for holomorphic self-maps of some fibered domains. \emph{Math. Res. Lett.}, \textbf{28}(3):697--706, 2021. MR4270270.
\bibitem[GM97]{Gen-Mig 1997} G. Gentili and S. Migliorini. A boundary rigidity problem for holomorphic mappings. \emph{Gen. Math.}, \textbf{5}:161--174, 1997. MR1723606.
\bibitem[Hua94]{Hua 1994} X. Huang. A non-degeneracy property of extremal mappings and iterates of holomorphic self-mappings. \emph{Ann. Scuola Norm. Sup. Pisa Cl. Sci. (4)}, \textbf{21}(3):399--419, 1994. MR1310634.
\bibitem[Hua95]{Hua 1995} X. Huang. A boundary rigidity problem for holomorphic mappings on some weakly pseudoconvex domains. \emph{Canad. J. Math.}, \textbf{47}(2):405--420, 1995. MR1335086.
\bibitem[HW22]{Hua-Wan 2022} X. Huang and X. Wang. Complex geodesics and complex Monge--Amp\`ere equations with boundary singularity. \emph{Math. Ann.}, \textbf{382}(3-4):1825--1864, 2022. MR4403236.
\bibitem[JP13]{Jar-Pfl 2013} M. Jarnicki and P. Pflug. \emph{Invariant Distances and Metrics in Complex Analysis}. Second extended edition. De Gruyter Expositions in Mathematics, vol.~9. Walter de Gruyter, Berlin, 2013. MR3114789.
\bibitem[JPZ93]{Jar-Pfl-Zei 1993} M. Jarnicki, P. Pflug, and R. Zeinstra. Geodesics for convex complex ellipsoids. \emph{Ann. Scuola Norm. Sup. Pisa Cl. Sci. (4)}, \textbf{20}(4):535--543, 1993. MR1267598.
\bibitem[KW13]{Kos-War 2013} \L. Kosi\'nski and T. Warszawski. Lempert theorem for strongly linearly convex domains. \emph{Ann. Polon. Math.}, \textbf{107}(2):167--216, 2013. MR3007944.
\bibitem[Lem81]{Lem 1981} L. Lempert. La m\'etrique de Kobayashi et la repr\'esentation des domaines sur la boule. \emph{Bull. Soc. Math. France}, \textbf{109}(4):427--474, 1981. MR0660145.
\bibitem[Lem82]{Lem 1982} L. Lempert. Holomorphic retracts and intrinsic metrics in convex domains. \emph{Anal. Math.}, \textbf{8}(4):257--261, 1982. MR0690838.
\bibitem[Lem84]{Lem 1984} L. Lempert. Intrinsic distances and holomorphic retracts. In \emph{Complex Analysis and Applications '81 (Varna, 1981)}, pages 341--364. Publ. House Bulgar. Acad. Sci., Sofia, 1984. MR0883254.
\bibitem[LQWZ26]{Lou-Qia-Wan-Zho 2026} C. Lou, J. Qiao, H. Wang, and Y. Zhong. A Burns--Krantz theorem on convex domains. Preprint, arXiv:2611.22516 [math.CV], 2026.
\bibitem[Mer93]{Mer 1993} P.~R. Mercer. Complex geodesics and iterates of holomorphic maps on convex domains in $\mathbb C^n$. \emph{Trans. Amer. Math. Soc.}, \textbf{338}(1):201--211, 1993. MR1123457.
\bibitem[PZ05]{Pfl-Zwo 2005} P. Pflug and W. Zwonek. Description of all complex geodesics in the symmetrized bidisc. \emph{Bull. London Math. Soc.}, \textbf{37}(4):575--584, 2005. MR2143737.
\bibitem[Ran95]{Ran 1995} T. Ransford. \emph{Potential Theory in the Complex Plane}. London Mathematical Society Student Texts, vol.~28. Cambridge University Press, Cambridge, 1995.
\bibitem[Ron26]{Ron 2026} F. Rong. A brief survey on the boundary rigidity. \emph{J. Geom. Anal.}, \textbf{36}(3):Paper No. 113, 15 pp., 2026. MR5032482.
\bibitem[RW83]{Roy-Won 1983} H.~L. Royden and P.-M. Wong. Carath\'eodory and Kobayashi metric on convex domains. Preprint, 1983.
\bibitem[Zim22]{Zim 2022} A. Zimmer. Two boundary rigidity results for holomorphic maps. \emph{Amer. J. Math.}, \textbf{144}(1):119--168, 2022. MR4367416.
\bibitem[Zwo25a]{Zwo 2025a} W. Zwonek. Burns--Krantz rigidity in non-smooth domains. \emph{Canad. Math. Bull.}, \textbf{68}(4):1037--1046, 2025. MR4992551.
\bibitem[Zwo25b]{Zwo 2025b} W. Zwonek. On left inverses in the Lempert theorem. Preprint, arXiv:2511.22301 [math.CV], 2025.
\end{thebibliography}
